\documentclass[11pt]{amsart} 
\usepackage{fullpage}

\usepackage{graphicx}

\usepackage{amsmath,amsthm,amssymb}
\usepackage{xcolor}
\usepackage{url}
\usepackage{ mathrsfs }

\newcommand     {\abs}[1]       {{\left\lvert{#1}\right\rvert}}

\newcommand{\std}{\textrm{std}}

\newcommand{\cond}{\textrm{cond}}

\newcommand{\legendre}[2]{\genfrac{(}{)}{}{}{#1}{#2}}

\newcommand{\Q}{\ensuremath{\mathbf{Q}}}

\newcommand{\cD}{\ensuremath{\mathcal{D}}}

\newcommand{\Z}{\ensuremath{\mathbf{Z}}}

\newcommand{\ep}{\epsilon}

\newcommand{\D}{\mathcal D_f}
\newcommand{\F}{\ensuremath{\mathcal{F}}}

\newcommand\be{\begin{equation}}
\newcommand\ee{\end{equation}}
\newcommand\bi{\begin{itemize}}
\newcommand\ei{\end{itemize}}
\newcommand\ben{\begin{enumerate}}
\newcommand\een{\end{enumerate}}

\newtheorem{lemma}{Lemma}
\newtheorem{theorem}[lemma]{Theorem}
\newtheorem{defn}[lemma]{Definition}

\newtheorem{corollary}[lemma]{Corollary}

\begin{document}

\title{Repulsion of zeros close to $s=1/2$ for L-functions }

\author[Coloma]{Nicolás Coloma}
\address{University of Colorado, Boulder, Colorado, USA}
\email{Nicolas.ColomaCarphio@colorado.edu}

\author[Espericueta]{Maria Espericueta Sandoval}
\address{Texas A\&M University--International, Laredo, Texas, USA}
\email{mariaespericuetasand@dusty.tamiu.edu}

\author[Lopez]{Erika Lopez}
\address{California State University--Channel Islands, Camarillo, California, USA}
\email{erika.lopez544@myci.csuci.edu}

\author[Ponce]{Francisco Ponce}
\address{North Carolina State University, Raleigh, North Carolina, USA}
\email{fmponcec@ncsu.edu}

\author[Rama]{Gustavo Rama}
\address{Universidad de la Rep\'ublica, Montevideo, Uruguay}
\email{grama@fing.edu.uy} 

\author[Ryan]{Nathan C. Ryan}
\address{Bucknell University, Lewisburg, Pennsylvania, USA}
\email{nathan.ryan@bucknell.edu}

\author[Vargas]{Alejandro Vargas-Altamirano}
\address{Bucknell University, Lewisburg, Pennsylvania, USA}
\email{aava001@bucknell.edu}

\begin{abstract}
In this paper we present results of several experiments in which we model the repulsion of low-lying zeros of L-functions using random matrix theory.  Previous work has typically focused on the twists of L-functions associated to elliptic curves and on families that can be modeled by unitary and orthogonal matrices.  We consider families of L-function of modular forms of weight greater than 2 and we consider families that can be modeled by symplectic matrices.  Additionally, we explore a model for low-lying zeros of twists that incorporates a discretization arising from the Kohnen--Zagier theorem.  Overall, our numeric evidence supports the expectation that the repulsion of zeros decreases as the conductor of the twist increases.  Surprisingly, though, it appears that using the discretization that arises from the Kohnen--Zagier theorem does not model the data better than if the discretization is not used for forms of weight 4 or above.
\end{abstract}

\maketitle

\section{Introduction}

The statistical properties of zeros and central values of L-functions have been extensively studied, computationally, heuristically and analytically.  A fruitful approach to studying these statistical properties has been to associate the statistics of an ensemble of matrices from a classical matrix group to the statistics of an L-function (or a collection of L-functions).  For instance, Montgomery \cite{montgomery} famously conjectured a formula for the pair correlation of the nontrivial zeros of the Riemann zeta function that is the same as the formula for the pair correlation for the eigenvalues of random matrices taken from either the Circular Unitary Ensemble or the Gaussian Unitary Ensemble if one takes the limit as the matrix size goes to infinity.  Odlyzko \cite{odlyzko} carried out massive computations of zeros of the Riemann zeta function to verify this conjecture.  Bogomolny and Keating \cite{bogomolnykeating1,bogomolnykeating2} provided heuristic evidence that not only the pair correlation functions agree, but all the $n$-point statistics do, as well.  Finally, Rudnick and Sarnak \cite{rudnicksarnak} proved that the zeros of the Riemann zeta function and the eigenvalues of this random matrix ensemble have the same $n$-point statistics in a restricted range.  

After studying the relationship between the zeros of the zeta function and the eigenvalues of random unitary matrices, analogous work was done with other L-functions and matrix groups.  In particular, inspired by the work of Katz and Sarnak \cite{katzsarnak}, collections of L-functions could be placed in families and those families could be associated to matrices from classical matrix groups.  Moreover, in this case, it appeared that the zeros of the L-functions had the same statistics as the eigenvalues of a randomly chosen matrix from the group, at least up to the leading term.  After the work of Katz and Sarnak there quickly appeared many examples of L-functions families behaving in a manner predicted by random matrix theory. Some of the families considered were: L-functions associated to holomorphic cusp forms (in either weight or level aspect); Dirichlet L-functions (either all or quadratic); and various twists or symmetric powers of L-functions. Analytic results about low-lying zeros were shown by Iwaniec, Luo, and Sarnak \cite{ils}, Rubinstein \cite{rubinstein}, \"Ozl\"uk and Snyder \cite{os}, and others. 

In this paper we focus on families of quadratic twists of L-functions of holomorphic modular forms; in particular, we want to understand the distribution of their low-lying zeros.  Miller~\cite[Figures 3 and 4]{miller06} observed that the first normalized zero above the central point of L-functions attached to rank-0 elliptic curves was repulsed from the central point. Due\~nez, Huynh, Keating, Miller, and Snaith recorded~\cite[Figure 4]{DHKMS} a similar repulsion in the family of even quadratic twists of the elliptic curve $E_{11}$.  Also in \cite{DHKMS}, they considered an ``excised'' model in which, because of a zero free region near $s=1/2$ guaranteed by theorems of Kohnen--Zagier \cite{KZ}, Baruch--Mao \cite{BaruchMao} and Mao \cite{Mao}, they only consider matrices with eigenvalues above a certain cutoff.  

In a recent preprint of Barrett and Miller \cite{bm}, similar analytic work is done for families of quadratic twists that do not correspond to the orthogonal group as families of twists of elliptic curves do.  In this paper, we experimentally study the repulsion of low-lying zeros for these families.  The scale of the experiment we carry out allow for the verification of conjectures relating random matrices and low-lying zeros of twists of L-functions corresponding to forms to which we cannot apply the Kohnen--Zagier theorem and to forms to which we can.  We also study the low-lying zeros of twists of higher weight forms which is of interest because the repulsion in weight 2 might be coming from there being infinitely many central values that vanish when we twist a weight 2 form; fewer but infinitely many when we twist a fom of weight 4, and finitely many (or none) when we twist a form of weight 6 or higher.  Our last collection of experiments is to extend the excised model to weights larger than 2 and compare how well the excised model and the non-excised model describe the repulsion of the low-lying zeros.

The paper is organized as follows.  In the next section, we give the necessary random matrix theory and L-function background and we summarize the main results in \cite{bm}.  In the subsequent section we present qualitative results of how well the different groups of random matrices model the repulsion of low-lying zeros for each of our families.  In the final section of the paper we discuss the excised model mentioned above and the computations we carried out to test the model.  We find that the excised model does not appear to do better than the standard model for weights bigger than 2.  See Figure~\ref{fig:cutoff_std}.  We also provide a possible explanation for why this might be the case.

\subsection*{Acknowledgments}  The work of Espericueta, Lopez, Ryan and Vargas was supported by the NREUP program of the Mathematical Association of America funded by the NSF Grant \#DMS-1950644.

\section{Background}

Let $f(q)=\sum_{n=1}^{\infty} a_n q^n \in S_k(M,\chi)$ be a classical newform of weight $k$, level $M$, character $\chi$ and let $\lambda_n=a_n/\sqrt{n}^{k-1}$.  In what follows we only consider normalized Hecke eigenforms (these are sometimes called primitive forms).  For $D>0$ a fundamental discriminant, let 
\[ 
L(f,s,\psi_D)=\sum_{n=1}^{\infty}  \psi_D(n)\frac{\lambda_n}{n^s} 
\]
be the L-series of $f$ (in the analytic normalization) twisted by the quadratic character $\psi_D$ associated to the real quadratic field $\Q(\sqrt{D})$; it will be useful to think of this character $\psi_D$ as the Kronecker symbol $\legendre{D}{\cdot}$.  Since we are restricting to fundamental discriminants $D$, the character $\psi_D$ is primitive.  The L-series has an analytic continuation $\Lambda(f,s,\psi_D)$ to the whole complex plane that satisfies the functional equation
\begin{equation}\label{eqn:FE}
\Lambda(f,s,\psi_D) = \epsilon_f\chi_f(D)\psi_D(-M)\overline\Lambda(f,1-s,\psi_D)
\end{equation}
for some complex number $\epsilon_f$ on the unit circle that only depends on $f$; here $\overline{\Lambda}(s) = \overline{\Lambda(\overline{s})}$.  The sign $\epsilon_f\chi_f(D)\psi_D(-M)$ will allow us to separate our L-functions into three different families; see Section~\ref{sec:families}.

The central values $L(f,1/2,\psi_D)$ can encode interesting arithmetic information about the form $f$, and a number of explicit investigations have been carried out examining the family of these values \cite{BSP,Gross,PT2,PT1,MRVT}.  These values play an important role on their own but they also provide connection between L-functions and random matrix theory.  The zeros of L-functions are particularly interesting because the location of their zeros is related to how many arithmetic objects of a certain size and kind there are.

\subsection{Modular forms}\label{sec:mfs}

In what follows, we will consider L-functions attached to three different kinds of modular forms.  Let $f(q) = \sum_{n=1}^\infty a_nq^n \in S_k(M, \chi)$ be a newform of weight $k$, odd level $M$, and character $\chi$. Then three cases emerge:  \begin{enumerate}
\item\label{mf:regular} $f$ could have principal character, 
\item\label{mf:nsd} $f$ could have non-trivial character $(f\neq \overline{f})$, or
\item\label{mf:sd} $f$ could have complex multiplication by its own non-trivial character $(f=\overline{f})$.  
\end{enumerate}
In the above list, according to \cite{bm}, forms of type~\ref{mf:regular} should have L-functions whose quadratic twists have zeros that are modeled by random matrices from the orthogonal group, forms of type~\ref{mf:nsd} should have L-functions whose zeros are modeled by matrices from the unitary group and forms of type~\ref{mf:sd} should have L-functions whose zeros are modeled by matrices from the symplectic group.  A newform is self-dual if its Fourier coefficients are real. A newform has complex multiplication (i.e., is CM) if there is a nontrivial Dirichlet character $\eta$ such that $\eta(p)a(p) = a (p)$ for all primes $p$ in a set of primes of density 1.  A form that is self-CM (as defined in \cite{bm}) is a form that is both self-dual and CM.  In order to make things more concrete, the particular modular forms we will consider are listed in Table~\ref{tbl:mfs}.

\begin{table}\small
\begin{tabular}{l|rr}
LMFDB label & Fourier expansion & Type \\\hline\hline
\texttt{11.2.a.a} & $f(q)=  q - 2q^{2} - q^{3} + O(q^{4})$ & $\chi$ principal\\
\texttt{7.4.a.a} & $f(q)=q - q^{2} - 2q^{3} + O(q^{4})$ & $\chi$ principal\\
\texttt{3.6.a.a} & $f(q)=q - 6q^{2} + 9q^{3} + O(q^{4})$ & $\chi$ principal\\
\texttt{3.8.a.a} & $f(q) = q + 6q^{2} - 27q^{3} + O(q^{4})$ & $\chi$ principal\\\hline\hline
%\texttt{3.10.a.b} & $f(q) = q + 18 q^{2} + 81 q^{3} + O(q^{4}) $ & $\chi$ principal \\
\texttt{13.2.e.a} & $f(q)= q + ( -1 - \zeta_{6} ) q^{2} + ( -2 + 2 \zeta_{6} ) q^{3} + O(q^{4})$ & $f\neq \overline{f}$\\\hline\hline
\texttt{7.3.b.a} & $f(q)= q - 3q^{2} + O(q^{4})$ & self-CM\\\hline\hline
\end{tabular}
\caption{Particular modular forms we will be using our experiments.  Here $\zeta_6$ is a particular $6$th root of unity.}\label{tbl:mfs}
\end{table}

\subsection{Admissible discriminants and families of L-functions}\label{sec:families}

Using the notation and terminology from above, we can now define the families of twists we will be computing and comparing to the predictions from random matrix theory.

\begin{defn}\label{def:D}
  Let $\mathcal D$ denote the set of fundamental discriminants.  Let $f\in S_k(M,\chi_f)$ be a newform, with $M$ an odd prime. If $f$ is self-CM, assume $\ep_f=+1$.  With these restrictions on $f$, let $\heartsuit\in\{\pm1\}$ and $1\leq\diamondsuit<M$ be integers, and put
  \be\label{eq:dgooddef}\D(X):=\left\{\begin{aligned}
      &\left\{D\in\cD:0<D\leq X\text{ and }\psi_D(M)\ep_f=+1\right\}
      &&\text{$\chi_f$ principal,} \\
      &\left\{D\in\cD:0<D\leq X\text{ and }\psi_D(M)=\heartsuit \right\}
      &&\text{$f$ self-CM,} \\
      &\left\{D\in\cD:0<D\leq X\text{ and }D\equiv\diamondsuit\mod M \right\}
      &&\text{$f\ne\overline f$.}\end{aligned}\right.\ee
\end{defn}
We now make precise our family $\F_f$.
\begin{defn}\label{def:F}
  Let \be\F_f(X):=\left\{L(f,s,\psi_D):D\in\D(X)\right\},\ee
  with $f\in S_k(M,\chi_f)$ as in Definition~\ref{def:D}.
\end{defn}
Then, if $\chi_f$ is principal, $\F_f$ is the family of quadratic twists of $L_f(s)$ where roughly half the central values of the twists vanish and the other half do not, depending on the parity of $\psi_D$. If $f$ is self-CM, $\F_f$ is the subfamily of the family of even quadratic twists, with an added condition on the $\psi_D(M)$ so that there are fewer vanishing central values among the twists.  If $f\ne\overline f$, then there is no notion of the parity of the functional equation of $L_f(s)$, and $\F_f$ is a subfamily of the family of quadratic twists of $L_f(s)$, with an added condition on the residue class of $D\mod M$.

\subsubsection{Justification of the congruence conditions}

In this section we study the signs of the functional equations of the three types of modular forms in order to motivate the definitions of the families in Definition~\ref{def:D}.  We also state and prove a nice theorem about the central values of twists that is related to the definition the unitary family.

We first consider the $f\neq \bar f$ case.  Let $f\in S_k(\Gamma_0(M),\chi)$, with $\chi$ primitive and nontrivial, and $L(f,s)$ be its L-function. The family of quadratic twists of this form is expected to have unitary symmetries and so we use the subscript $U$.  Define, for $0<\diamondsuit<M$,
\[
\cD_U(f,\diamondsuit) := \{\psi_D = \legendre{D}{\cdot}: D \text{ fundamental and }  0 < D = \cond(\psi_D)\equiv \diamondsuit\pmod{M}\}
\]
and the associated family of twisted L-functions to be
\[
\F_U(f,\diamondsuit) := \{L(f,s,\psi): \psi\in \cD_U(f,\diamondsuit)\}.
\]
The completion of an L-function in $\F(f,\diamondsuit)$ satisfies the functional equation \eqref{eqn:FE} that we recall here:
%\begin{equation}\label{eqn:symplectic-twisted-FE}
\[
\Lambda(f, s, \psi) = \chi(D)\psi (-M)\epsilon_f \overline{\Lambda}(f, 1-s,\psi)
%\end{equation}
\]
with $\epsilon_f$ a complex number on the unit circle that only depends on $f$.  The fact that $\chi(D)\psi_D(-M)\epsilon_f$ is essentially an arbitrary complex number on the unit circle suggests that the flatness of the second plot in Figure~\ref{fig:densities} is reasonable.

In order to justify the need to separate discriminants according to congruence classes in the unitary case, we show:
\begin{theorem}  Let $f\in S_k(\Gamma_0(M),\chi)$ with $\chi$ nontrivial and let $\diamondsuit$ be an integer so that $1\leq \diamondsuit < M$.  Then all the central values $L(f,1/2,\psi)$ for $L(f,s,\psi)\in \F_U(f,\diamondsuit)$, lie on a line through the origin.
\end{theorem}

\begin{proof}

  Let $D>0$ be a fundamental discriminant and $\psi_D \in \cD(f,\diamondsuit)$ be arbitrary; i.e., suppose $\psi_D$ is a Kronecker symbol $\legendre{D}{\cdot}$ where $\cond(\psi)=D\equiv \diamondsuit\pmod{M}$. Now consider $L(f,s,\psi) \in \F(f,\diamondsuit)$. Then, by \eqref{eqn:FE} we have
  \[
 \Lambda(f, 1/2, \psi_D) = \chi(D)\psi_D (-M)\epsilon_f \overline{\Lambda}(f,1/2, \psi).
  \]
  Thus,
  \[
  \frac{\Lambda(f,1/2, \psi_D)}{\overline{\Lambda}(f, 1/2, \psi_D)} = \chi(D)\psi_D (-M)\epsilon_f.
  \]
  Since the character $\psi_D$ is real and the coefficients of $f$ are complex, we know $\overline{\Lambda}(f,1/2,\psi)=\overline{\Lambda(f,1/2,\psi)}$ and thus
  \[
 \frac{\Lambda(f,1/2, \psi_D)}{\overline{\Lambda(f, 1/2, \psi_D)}} = \chi(D)\psi_D (-M)\epsilon_f.
  \]
So the arguments of $\Lambda(f,1/2,\psi)$ and $\sqrt{\chi(D)\psi (-M)\epsilon_f}$ differ by an integer multiple of $\pi$. 
    
  Now $\chi$ is periodic with period $M$; thus, since $D\equiv \diamondsuit\pmod{M}$, we know $\chi(D) = \chi(\diamondsuit)$.  In particular, it only depends on the congruence class of $D$ mod $M$.
  Next,
  \[
  \psi_D(-M) = \legendre{D}{-M}=\legendre{D}{-1}\legendre{D}{M}=\legendre{D}{M}=\legendre{\diamondsuit}{M}=:s_\diamondsuit
  \]
  since $M$ is odd.  In particular it only depends on the congruence class of $D$ mod $M$.

%    since $\cond(\psi) \equiv d$ mod $N$ and $p_{i}$ is a prime factor of $N$. Thus $\psi(N)$ is only dependent on the congruence class of $\cond(\psi)$ mod $N$, and we can replace $\psi(N)$ with some $s_{d}$ which only depends on $d$.

%Finally, observe
%\[
%\cond(\psi)=|\tau(\psi)|^2 = \tau(\psi)\overline{\tau(\psi)} = \tau(\psi)\psi(-1)\tau(\bar\psi)=\tau(\psi)\psi(-1)\tau(\psi)=\psi(-1)\tau(\psi)^2
%\]
%since $\psi$ is quadratic and, therefore, real.  Moreover, since all our twists are by characters corresponding to imaginary quadratic fields, we know that $\psi(-1)= \legendre{D}{-1}=-1$,

%  Next note that
%  \[
%  \frac{\tau(\psi)^{2}}{\cond(\psi)} = \frac{\tau(\psi){\tau(\bar{\psi})}}{\cond(\psi)} = \psi(-1)\frac{\tau(\psi)\overline{\tau(\psi)}}{\cond(\psi)},
%  \]
%  but we know that $\tau(\chi)\overline{\tau(\chi)} = c$, and since $c$ is positive, we have $\chi(-1) = 1$, thus $\frac{\tau(\chi)^{2}}{c} = 1$.

Taking this all together we are left with
\[
\sqrt{\chi(-D)\psi_D (-M)\epsilon_f}= \sqrt{\chi(\diamondsuit)s_\diamondsuit\epsilon_f}.
\]
But the right-hand side of the equation is constant for all L-functions in $\F_U(f,\diamondsuit)$. And since the argument of the central value of any L-function in $\F_U(f,\diamondsuit)$ will differ from the argument of $\sqrt{\chi(\diamondsuit)s_{\diamondsuit}\epsilon_f}$ only by an integer multiple of $\pi$, we can conclude that the arguments of the central values of all L-functions in $\F_U(f,\diamondsuit)$ differ only by a multiple of $\pi$ which means the central values must be on a line through the origin in the complex plane, with slope $\tan(\arg(\sqrt{\chi(\diamondsuit)s_{\diamondsuit}\epsilon_f}))$.
\end{proof}

%Now, let $f\in S_k(\Gamma_0(M),\chi)$, with $\chi$ principal, and $L(f,s)$ be its L-function. The family of quadratic twists of this form is expected to have orthogonal symmetries and so we use the subscript $O$.  Define, for $\diamondsuit>0$,
%\[
%\cQ_O(f,\diamondsuit) := \{\psi \text{ primitive quadratic character}: \cond(\psi)\equiv \diamondsuit\pmod{M}\}
%\]
%and the associated family of twisted L-functions to be
%\[
%\cL_O(f,\diamondsuit) := \{L(f,s,\psi): \psi\in \cQ_O(f,\diamondsuit)\}.
%\]

So, in the unitary case, the choice of a congruence class mod $M$ is reasonable.  In the other two cases, we know that the central values are real because $\Lambda(f,s,\psi_D)=\overline{\Lambda(f,s,\psi_D)}$ in those cases since $f$ has real coefficients and $\psi_D$ only takes on real values.  So, we can state the following corollary of perhaps independent interest:
\begin{corollary}
Let $f$ belong to one of the three classes of modular forms listed in Section~\ref{sec:mfs}.  Then the central values $L(f,1/2,\psi_D)$ for $D$ as in Definition~\ref{def:D} all lie on a line through the origin.
\end{corollary}

We still need to justify the other two conditions in Definition~\ref{def:D}.  In the case when $\chi$ is principal, we know that the sign of the functional equation becomes
\[
\chi(-D)\psi_D(-M)\epsilon_f = \psi_D(-1)\psi_D(M)\epsilon_f
\]
and so if $\psi_D(M)\epsilon_f=1$, then the central value vanishes whenever $\psi_D$ is odd.  This explains why the peak in the middle of the first plot in Figure~\ref{fig:densities} is reasonable.

When $f$ is self-CM, we restrict to $\epsilon_f=1$ in Definition~\ref{def:D}.  In \cite[Lemma 3]{bm} it is shown that $\epsilon_f \chi_f(D)\psi_D(-M)=\epsilon_f$ when $f$ is self-CM.  In other words, there are no forced vanishings at $s=1/2$.  This explains why the dip in the middle of the third plot in Figure~\ref{fig:densities} is reasonable.

\subsubsection{Theorem of Kohnen--Zagier}
The collection $L(f,1/2,\psi_D)$ of central values with varying discriminant $D$ for self-dual $f$ that are not CM can be computed efficiently via a theorem of Kohnen--Zagier \cite{KZ} and generalizations to higher level by Baruch-Mao \cite{BaruchMao} and Mao \cite{Mao}.  This theorem asserts that the central values are related to the Fourier coefficients of a certain half-integer weight modular form.  Concretely, for a fundamental discriminant $D>0$ coprime to $M$, we have
\begin{equation} \label{eqn:LfchiD}
L(f,1/2,\psi_D) = \kappa_f  \frac{c_{D}(g)^2}{\sqrt{D}^{k-1}}
\end{equation}
where the (nonzero) constant $\kappa_f$ is independent of $D$ and the integer $c_{D}(g)$ is the $D$th coefficient of a modular form $g$ of weight $(k+1)/2$ related to $f$ via the Shimura correspondence.  We point out an interesting consequence of this theorem:  if $L(f, 1/2, \psi_D) < \kappa_f  \frac{1}{\sqrt{D}^{k-1}}$, then $L(f, 1/2, \psi_D)=0$ because the coefficients of $g$ are integral.  This gives the set of central values of twists a discretization, at least for those $f$ to which the above theorem can be applied.  This will be used in Section~\ref{sec:excised} when we examine the excised model mentioned above.

%Computing central values using \eqref{eqn:LfchiD} has the advantage that the description of $g$ as a linear combination of theta series permits the rapid computation of a large number of coefficients: for example, Hart--Tornar\'ia--Watkins \cite{hart2010congruent} compute hundreds of billions of twists of the congruent number elliptic curve (using FFT methods).  Random matrix theory has proved useful in refining conjectures related to the low-lying zeros of L-functions \cite{KeatingSnaith2,KeatingSnaith1}.  In work of Conrey--Keating--Rubinstein--Snaith \cite{CKRS}, the following basic question was considered.  Let $f \in S_k(M)$ be a newform with rational integer coefficients.  For how many fundamental discriminants $D$ with $|D|\leq X$ does the twisted L-function  $L (f,s,\psi_D)$ vanish at the center of the critical strip?  In a collection of papers \cite{MRVT,PT,TR}, a number of variants of this problem were considered: the weight of $f$ was allowed to vary, the level of $f$ was allowed to be composite, and so on.  We extend the excised model to weights beyond 2.

\subsubsection{Repulsion of the lowest zero}

In \cite{miller06} it was observed that the first zeros of elliptic curve L-functions in certain families exhibited a repulsion of the first zero above the real line.  Since the eigenvalues of matrices in $SO(2N)$ do not exhibit this same repulsion, this runs counter to the expectation that the statistics of eigenvalues of random matrices should be similar to the statistics of zeros of L-functions in a family.  The expectation is believed to hold in the limit as the conductor of the L-functions tend to infinity (and the size of the matrices do, too), but there can be discrepancies for ``finite'' conductors.

It has been proposed \cite{DHKMS} to restrict the choices of matrices from $SO(2N)$ so that their characteristic polynomials evaluated at 1 cannot be too small without being zero in much the same way that the Kohnen--Zagier Theorem and its generalizations described above says that an L-function's central value cannot be too small without being zero.  This is the so-called ``excised'' model mentioned above.  In \cite{DHKMS} some computational evidence was provided that this model fits the L-function data better for twists of elliptic curve L-functions and in \cite{marshall}, some theoretical evidence was given that, for families of twists of elliptic curve L-functions, there really is repulsion of the first zero above the real line.

In what follows we examine the repulsion of families modeled by $U(N)$ and $USp(2N)$ and also the repulsion of families modeled by $SO(N)$ but where the families are twists of modular form L-functions for weight $>2$.  Finally, we consider the effectiveness of the excised model for orthogonal families of twists of modular form L-functions for weights $>2$.  In the next section we introduce these groups of random matrices.

\subsection{Random matrices}

As shown in \cite{KeatingSnaith2,KeatingSnaith1,rubinstein,montgomery} and elsewhere, the local statistical properties of the Riemann zeta function and other L-functions can be modeled by the characteristic polynomials of Haar distributed random matrices.  There are three groups that we will use in what follows and that we have already referred to above:  the unitary group $U(N)$, the even special orthogonal group $SO(2N)$ and the unitary symplectic group $USp(2N)$.  These groups are made into probability spaces by using each group's Haar measure as the space's distribution.  

\begin{figure}
\includegraphics[width=.5\textwidth]{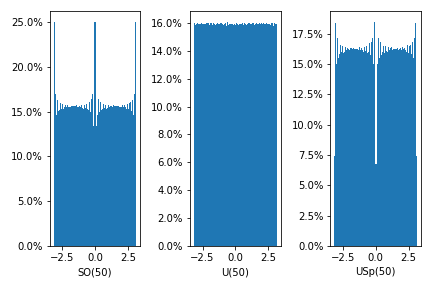}
\caption{Density plots of 100,000 $50\times 50$ matrices in (from the left) the special orthogonal group, the unitary and the unitary symplectic group.}\label{fig:densities}
\end{figure}

In order to carry out our experiments, we need to calculate large samples of random matrices from each of these groups.  We do this following \cite{mezzadri} and our implementation of the algorithm described there is available at \cite{code}.  See Figure~\ref{fig:densities} for density plots of $100,000$ $50\times 50$ matrices from each group; these distributions match the expected distributions and also match the vanishings of central values within families of twists as discussed heuristically in Section~\ref{sec:families}.

\subsubsection{Matrix sizes and cutoffs}

In order to compare the distributions of eigenvalues of random matrices and the low-lying zeros of L-functions, we have to determine the size of the matrices we use.  The standard approach to finding the matrix size is to choose the matrix size $N_\std$ so that the mean densities of eigenvalues are equal to mean density of zeros.  In particular, this means that for discriminants around $X$ we have $N_\std = \log\left ( \frac{\sqrt{3}X}{2\pi e}\right)$.  We point out that when we are modeling families using $SO(2N)$ or $USp(2N)$, we double $N_\std$.
%In this paper we study two different approaches to computing the matrix size and for each approach, we also determine the corresponding cutoff for the excised model.  See \cite{bm} for more details but we present the necessary results here.

%When we model twists that correspond to the orthogonal group, we also use an effective matrix size that tries to capture the lower order terms in these densities.  In particular, there is a constant $a_1$ (defined below), such that the $N_\eff = \tfrac{N_\std}{2a_1}$.

For the orthogonal group, we also analyze a model (first described in \cite{DHKMS}) of low-lying zeros of twists of L-functions that incorporates the discretization that comes from the Kohnen--Zagier theorem.  In particular, for a modular form of weight $k$, the zeros are discretized by $1/D^{(k-1)/2}$ and so we exclude from our random matrices those whose value at 1 is of the scale $\exp((1-k)N_\std/2)$; that is, we want matrices $A$ in $SO(2N)$ whose characteristic polynomials $\Lambda_A(z,N)$ satisfy
\[
\abs{\Lambda_A(1,N)} \geq c_\std\cdot \exp((1-k)N_\std/2)
\]
for some constant $c_\std$.  %The effective cutoff $c_\eff$ is defined similarly except we replace the $\std$ subscripts with $\eff$.  

In \cite{DHKMS} and \cite{bm} formulas are given for these cutoffs, but we follow the numerical method to estimate the cutoffs described in \cite{DHKMS}.  In this approach they try several values of $c_\std$ and empirically measure the distance between the cumulative distributions of zeros (these are independent of the choice of cutoff) and eigenvalues (these will be greater than or equal to the cutoff) by numerically approximating the area between them at several specified points.  The value of $c_\std$ for which this is smallest, is the value we use.  In \cite{DHKMS}, the value of $c_\std$ for quadratic twists of the L-function associated to the modular form \texttt{11.2.a.a} was computed to be $\approx 2.188$ and this was shown to agree with the formula for $c_\std$ in the same paper.  We follow their approach and get a different value of $c_\std$\footnote{Their numerical value of $c_\std$ might be wrong because the mean value of the first zero that they report is incorrect as verified by our code, by Rubinstein's lcalc \cite{lcalc} and PARI/GP \cite{pari}}.

\subsection{Computing the first few zeros of L-functions}

Our ultimate goal is to understand the distribution of zeros.  To calculate them, we use the implementation in PARI/GP \cite{pari} that is described in \cite{paribook}.  Roughly speaking, a naive search is done for zeros of the real-valued Hardy Z-function along the critical line $s=1/2$. The computations here are limited due to the complexity of calculating zeros of L-functions of large conductor.  In each case, we compute the first few zeros of twists of the L-functions of the orthogonal forms in Table~\ref{tbl:mfs} up to discriminant 1,000,000 and the symplectic and unitary forms up to discriminant 40,000. Our code and data are available at \cite{code}.

\section{The standard model}\label{sec:standard}

In this section we describe the results of our experiments using the standard random matrix model.  In particular, we show that the distribution of the first zero for each family of modular forms is roughly the same as the distribution of the argument of the first eigenvalue of a sample of random matrices from the corresponding matrix group.  

\subsection{Distributions}  
\begin{figure}
\includegraphics[width=.3\textwidth]{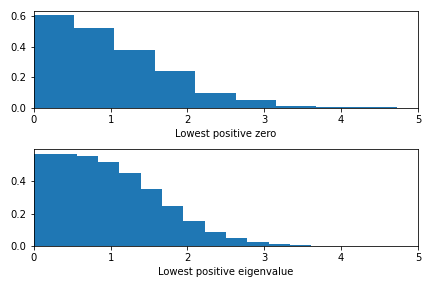}\hfill
\includegraphics[width=.3\textwidth]{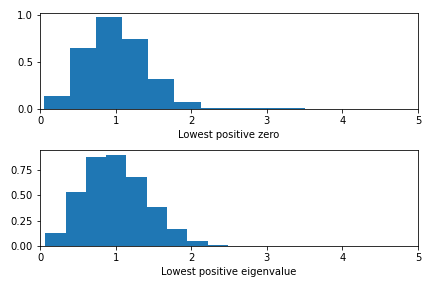}\hfill
\includegraphics[width=.3\textwidth]{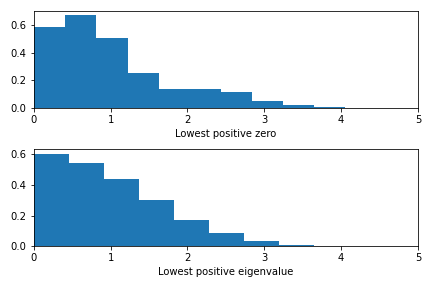}\\
\caption{Distributions of lowest zeros of admissible twists of \texttt{3.8.a.a} and lowest mean eigenvalues from $SO(2N)$ (left), distributions of lowest zeros of admissible twists of \texttt{13.2.e.a} and lowest mean eigenvalues from $USp(2N)$ (center), and distributions of lowest zeros of admissible twists of \texttt{7.3.b.a} and lowest mean eigenvalues from $U(N)$ (right). The data have been normalized to have a mean of one.}
\label{fig:dists-families}
\end{figure}

\begin{figure}
\hfill\includegraphics[width=.3\textwidth]{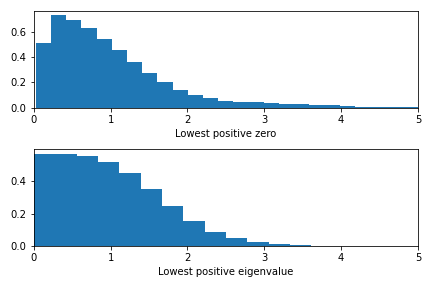}\hfill
\includegraphics[width=.3\textwidth]{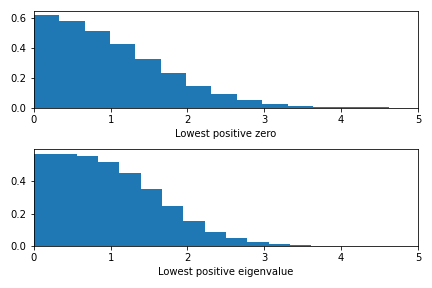}\hfill
\includegraphics[width=.3\textwidth]{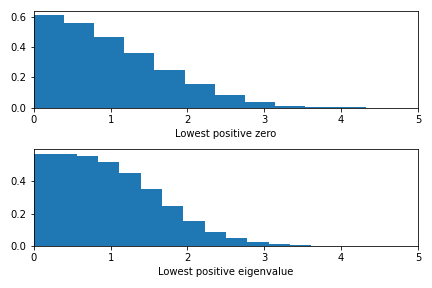}\hfill\\
\includegraphics[width=.3\textwidth]{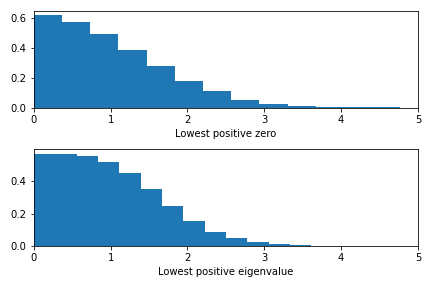}
\caption{Distributions of lowest zeros of admissible twists of forms in Table~\ref{tbl:mfs}.  From left to right the top plots in the first row are distributions of lowest zeros for admissible quadratic twists of the forms \texttt{11.2.a.a}, \texttt{7.4.a.a}, \texttt{3.6.a.a} and the top plots in the second row are distributions of lowest zeros for admissible quadratic twists of the form \texttt{3.8.a.a}. %and \texttt{3.10.a.b}.
  The bottom plots are the lowest eigenvalues from $SO(2N)$.  Both plots have been normalized to have means of 1. } 
\label{fig:dists-weights}
\end{figure}

Recall that the three types of families that we are considering are those that can be modeled by $SO(2N)$, $U(N)$ and $USp(2N)$ and in each case we have zeros for twists up to discriminant 40,000.  In this first experiment, we compare the distributions of the eigenvalues and zeros, normalized so that they both have means of 1 and observe that the shape of each pair of plots in Figure~\ref{fig:dists-families} and Figure~\ref{fig:dists-weights} are similar.  In our computations of the eigenvalues, we calculated $N_\std$ for each discriminant $D$ and then found the mean lowest eigenvalue over a sample of $10,000$ $N_\std \times N_\std$ matrices.

\subsection{Repulsion}

\begin{figure}
\includegraphics[width=.3\textwidth]{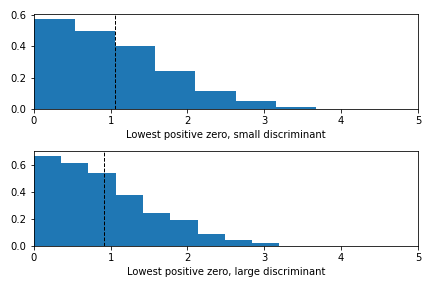}\hfill
\includegraphics[width=.3\textwidth]{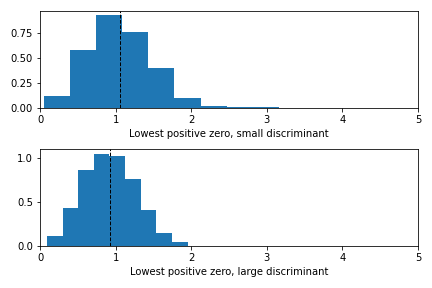}\hfill
\includegraphics[width=.3\textwidth]{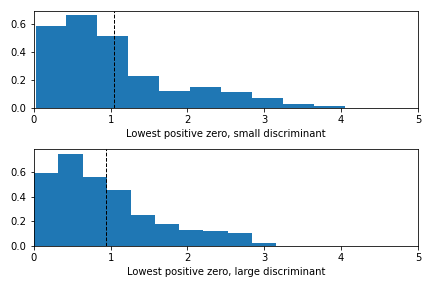}\\
\caption{Distributions of lowest zeros of admissible twists of \texttt{3.8.a.a} separated into those of small and large conductor (left), distributions of lowest zeros of admissible twists of \texttt{13.2.e.a} separated into those of small and large conductor (center), and distributions of lowest zeros of admissible twists of \texttt{7.3.b.a} separated into those of small and large conductor (right).  The dashed vertical lines in each graph are the means of the data; and, again, before splitting into small and large conductors, the data were normalized to have a mean of one.}\label{fig:repulsion}
\end{figure}

\begin{figure}
\includegraphics[width=.3\textwidth]{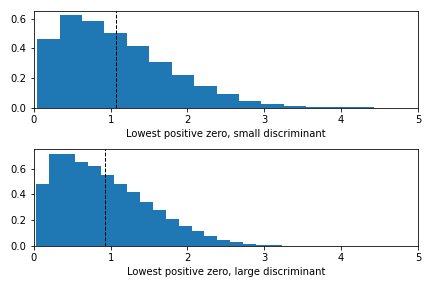}\hfill
\includegraphics[width=.3\textwidth]{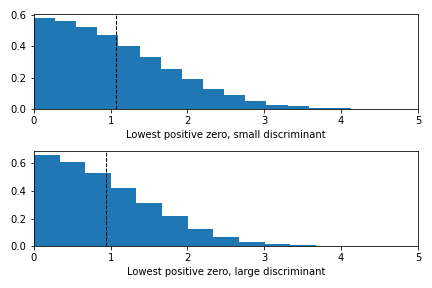}\hfill
\includegraphics[width=.3\textwidth]{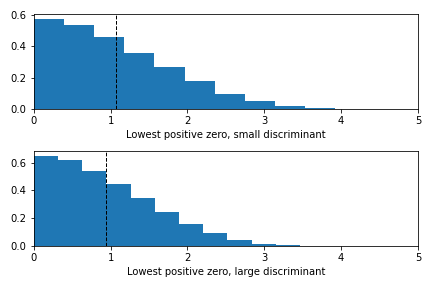}\\
\includegraphics[width=.3\textwidth]{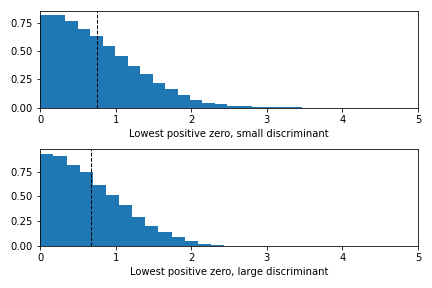}
\caption{Distributions of lowest zeros of admissible twists of \texttt{11.2.a.a} separated into those of small and large conductor (top row, left), distributions of lowest zeros of admissible twists of \texttt{7.4.a.a} separated into those of small and large conductor (top row, center), distributions of lowest zeros of admissible twists of \texttt{3.6.a.a} separated into those of small and large conductor (top row, right), distributions of lowest zeros of admissible twists of \texttt{3.6.8.a} separated into those of small and large conductor (bottom row, left)
  %and distributions of lowest zeros of admissible twists of \texttt{3.10.a.a} separated into those of small and large conductor (bottom row, right).
  The dashed vertical lines in each graph are the means of the data; and, again, before splitting into small and large conductors, the data were normalized to have a mean of one.}\label{fig:repulsion-weights}
\end{figure}

We verify the expectation that the average repulsion is less for larger conductors than it is for smaller conductors; see Figure~\ref{fig:repulsion} and Figure~\ref{fig:repulsion-weights}. For each of these plots we broke the set of admissible twists, ordered by discriminant, in half and called the first half as being of ``small'' conductor and the second half as being of ``large'' conductor.  According to the philosophy of the correspondence between zeros of L-functions and random matrix theory, as the discriminant of the twist goes to infinity, the repulsion goes to zero since this corresponds to the matrix size tending to infinity, and hence the smallest eigenvalue tends to 1, corresponding to an angle of 0.  In Figure~\ref{fig:repulsion} and Figure~\ref{fig:repulsion-weights} we observe this phenomenon: as predicted by random matrix theory the repulsion (as measured by the mean of the lowest zeros in each group) decreases with larger discriminants.

\section{The excised model}\label{sec:excised}

In this experiment we study how well the excised model does compared to the standard model that does not take into account the discretization from the Kohnen--Zagier theorem.  In order to do this we first need to find the cutoffs for the forms we are considering; because of the difficulty in computing central values of forms not modeled by the orthogonal group, we are limiting our attention to those forms in Table~\ref{tbl:mfs} for which $\chi$ is principal.

\subsection{Computing central values}

We need to numerically compute the distribution of the central values of the families of quadratic twists described in Table~\ref{tbl:mfs}.  We now describe how we do this in the cases when the family of quadratic twists has central values modeled by orthogonal matrices.  In order to calculate these central values for modular forms with principal character we use standard approaches but carry them out for a wider range of weights and to higher discriminant bounds.  Our data is available at \cite{code}. 

\subsubsection{Weight 2}  For the modular form $f\in S_2(11)$ with label \texttt{11.2.a.a}, we follow the method described in \cite{MRVT} to compute central values $L(f,1/2,\psi_D)$, using Brandt matrices.

\subsubsection{Weight 4}
The first extensive computations of the modular form $f\in S_4(7)$ with label \texttt{7.4.a.a} and its Shimura lift were carried out in \cite{TR}.   We proceed in a similar way (the details are slightly different to be consistent with our other computations) and compute its Fourier expansion as
\[f = \frac{1}{4}\sum_{(a, b, c, d)\in\Z^4}(2a^2+2ab-3b^2)q^{Q_7(a, b, c, d)}\]
where $Q_7(a, b, c, d) = a^2 + ab +2b^2 + c^2 + bc + 2d^2 = Q_7'(a, b) + Q_7'(c, d)$, with $Q_7'(x, y) = x^2 + xy + 2y^2$.
Then
\[f = \frac{1}{4}\left(\sum_{(a, b)\in\Z^2}(2a^2+2ab-3b^2)q^{Q'_7(a, b)}\right)\left(\sum_{(b, c)\in\Z^2}q^{Q'_7(b, c)}\right),\]
and we can compute the Fourier coefficients of $f$ in linear time.

To compute the half integral modular form $g_+$ associated to $f$ via the Shimura lift, we define
\[
w_{11}(x, y, z) = \left\{
\begin{array}{ll}
0&\text{if }11\nmid Q(x, y, z)\\
\left(\frac{-2x+z}{11}\right)&\text{if }2x\not\equiv z\pmod{11}\\
\left(\frac{x}{11}\right)&\text{otherwise}
\end{array}
\right.,
\]
and
\[g_+ = \frac{1}{4}\sum_{(x, y, z)\in\Z^3}xw_{11}(x, y, z)q^{Q(x, y, z)/11}=\sum_{n=1}^\infty c_+(n)q^n,\]
where $Q(x, y, z) = 4x^2+4xy+8y^2+7z^2$.
As before, we can compute this in linear time, this is because $xw_{11}(x, y, z)$ only depends on the variables $x$ and $z$,
and $Q(x, y, z) = Q'(x, y) + 7z^2$, where $Q'(x, y) = 4x^2+4xy+8y^2$.
So,
we define
\[
\sum_{n=1}^\infty a_n q^n = \frac{1}{4}\left(\sum_{(x, y)\in\Z^2}xw_{11}(x, y, 0)q^{Q'(x, y)}\right)\left(\sum_{z\in\Z}q^{7z^2}\right),
\]
so that
$c_+(n) = a_{11n}$.

\subsubsection{Weight 6}
For the modular form $f\in S_6(3)$ with label \texttt{3.6.a.a} we can compute its Fourier expansion as
\[f = \frac{1}{6}\sum_{(a, b, c, d)\in\Z^4}P(a, b, c, d)q^{Q_3(a, b, c, d)}\]
where $Q_3(a, b, c, d) = a^2 - ab + b^2 + c^2 - cd + d^2 = Q'_3(a, b) + Q'_3(a, b)$, $Q'_3(x, y) = x^2 - xy + y^2$, 
and $P(a, b, c, d) = a^4 - 2a^2c^2 - 2a^3b + 4ac^2b + 3a^2b^2 - 4b^2c^2 - 2ab^3 + b^4 - 2abcd + 4cb^2d - 2b^2d^2$.
So
\begin{align*}
f &= \\
\frac{1}{6}&\left(\sum_{(a, b)\in\Z^2}(a^4-2a^3b+3a^2b^2-2ab^3+b^4)q^{Q'_3(a, b)}\right)\left(\sum_{(c, d)\in\Z^2}q^{Q'_3(b, c)}\right)\\
+&\left(\sum_{(a, b)\in\Z^2}(-2a^2+4ab-4b^2)q^{Q'_3(a, b)}\right)\left(\sum_{(c, d)\in\Z^2}c^2q^{Q'_3(c, d)}\right)\\
-&\left(\sum_{(a, b)\in\Z^2}ab\,q^{Q'_3(a, b)}\right)\left(\sum_{(c, d)\in\Z^2}cd\,q^{Q'_3(c, d)}\right)\\
-&\left(\sum_{(a, b)\in\Z^2}2b^2\,q^{Q'_3(a, b)}\right)\left(\sum_{(c, d)\in\Z^2}d^2\,q^{Q'_3(b, c)}\right)\\
=&\left(\sum_{(a, b)\in\Z^2}(a^4-2a^3b+3a^2b^2-2ab^3+b^4)q^{Q'_3(a, b)}\right)\left(\sum_{(c, d)\in\Z^2}q^{Q'_3(b, c)}\right)\\
-&2\left(2\left(\sum_{(a, b)\in\Z^2}a^2q^{Q'_3(a, b)}\right)-\left(\sum_{(c, d)\in\Z^2}cd\,q^{Q'_3(c, d)}\right)\right)^2
\end{align*}

To compute the half integral modular form $g_+$ associated to $f$ we define
$Q_3(x, y, z) = 4x^2+4xy+4y^2 + 3z^2$,
\[
w_{7}(x, y, z) = \left\{
\begin{array}{ll}
0&\text{if }7\nmid Q_3(x, y, z)\\
\left(\frac{4x+5y}{7}\right)&\text{if }4x\not\equiv 5y\pmod{7}\\
\left(\frac{2x}{7}\right)&\text{otherwise}
\end{array}
\right.,
\]
\[
w_3(x, y, z) = \left(\frac{2x+y}{3}\right),
\]
and
\[
P(x, y, z) = 2x^2+2xy+2y^2-3z^2.
\]
Then
\[
g_+ = \frac{1}{6}\sum_{(x, y, z)\in\Z^3}P(x, y, z)w_3(x, y, z)w_7(x, y, z)q^{Q_3(x, y, z)}.
\]
Because $w_3$ and $w_7$ only depends in the variables $x, y$ we can compute the coefficients
of $g_+$ in linear time as before.

\subsubsection{Weight 8}
For the modular form $f\in S_8(3)$ with label \texttt{3.8.a.a} we can compute its Fourier expansion as
\[f = \frac{1}{6}\sum_{(a, b, c, d)\in\Z^4}P(a, b, c, d)q^{Q_3(a, b, c, d)}\]
with $Q_3$ as before, and
$P(a, b, c, d) = 2a^{6} - 6a^{5} b - 15a^{4} b^{2} + 40a^{3} b^{3} - 15a^{2} b^{4} - 6a b^{5} + 2b^{6} + 2c^{6} - 6c^{5} d - 15c^{4} d^{2} + 40c^{3} d^{3} - 15c^{2} d^{4} - 6c d^{5} + 2d^{6}
= P_1(a, b) + P_1(c, d)$.
So,
\[
f = \frac{1}{6}\left(\sum_{(a, b) \in\Z^2}P_1(a, b)q^{Q'_3(a, b)}
\right)
\sum_{(c, d)\in\Z^2}q^{Q'_3(c, d)}
.\]

To compute the half integral modular form $g_+$ associated to $f$ we define
$P(x, y, z) = 2x^3 + 3x^2y - 3xy^2 - 2y^3$,
and we have
\[
g_+ = \frac{1}{6}\sum_{(x, y, z)\in\Z^3}P(x, y, z)w_7(x, y, z)q^{Q_3(x, y, z)},
\]
and as before we can compute the coefficients of $g_+$ in linear time.

\subsection{Computing $c_\std$}

\begin{figure}
\includegraphics[width=.4\textwidth]{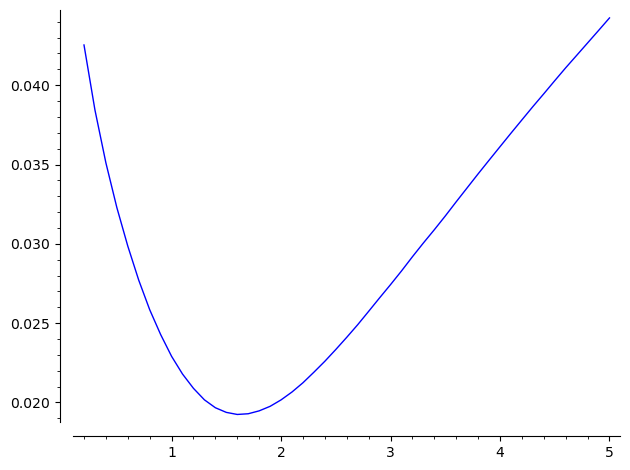}
\includegraphics[width=.4\textwidth]{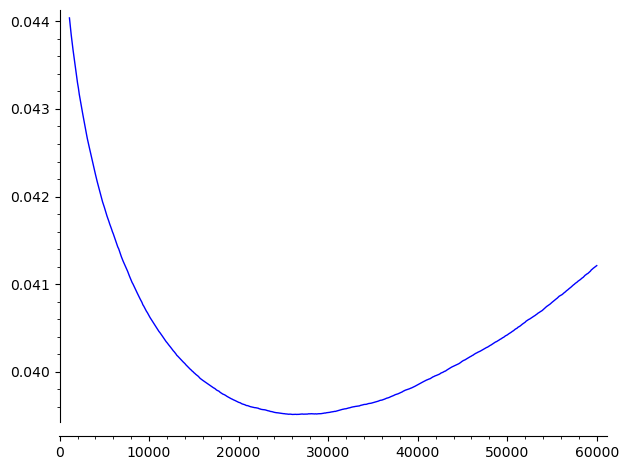}\\
\includegraphics[width=.4\textwidth]{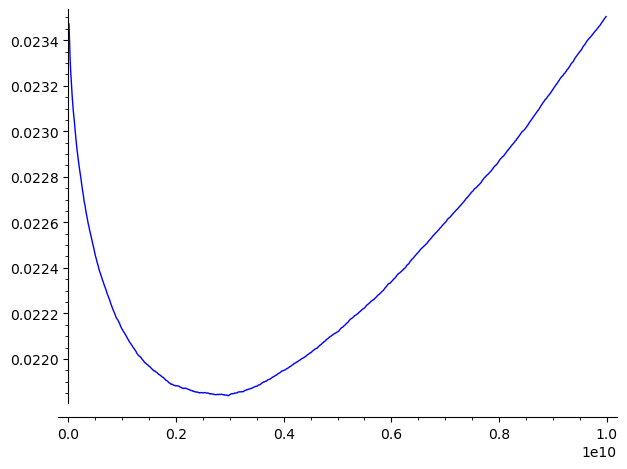}
\includegraphics[width=.4\textwidth]{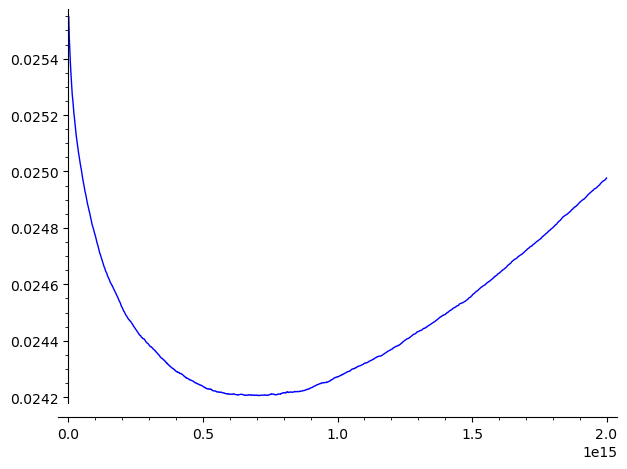}\\
\caption{Each of these four plots shows the difference, for various candidate values of the cutoff $c_\std$, between the cumulative distributions of values of $L(f,1/2,\psi_D)$ and of values of $\Lambda_A(1,N)$.  From left to right the forms whose central values $f$ are being calculated are \texttt{11.2.a.a}, \texttt{7.4.a.a}, \texttt{3.6.a.a}, and \texttt{3.8.a.a}.  In all cases, we generated random elements from $SO(24)$.}\label{fig:c_std}
\end{figure}

The formulas in \cite{bm} for $c_\std$ involves an asymptotic for the cardinality of the set of vanishings of the twists at $s=1/2$ and because in weights 6 and above there are expected to be no or finitely many such vanishings, we decided to proceed with a numerical approach.  As described above, we numerically approximate the value of $c_\std$ by comparing, for various candidates of $c_\std$, the cumulative distributions of central values cut off at the candidate value of $c_\std$ and the cumulative ditribution of evaluations of characteristic polynomials at 1 cut off at the candidate value of $c_\std$.  We make a plot of these differences for each candidate value of $c_\std$ and find the minimum on each plot.  As the weight increases, the computations of the central values get harder and so the plots are less smooth for larger weight.  See Figure~\ref{fig:c_std}.  

The value of $c_\std$ for twists of \texttt{11.2.a.a} was roughly $1.6$, the value of $c_\std$ for twists of \texttt{7.4.a.a} was roughly $2.6 \times 10^4$, the value of $c_\std$ for twists of \texttt{3.6.a.a} was roughly $2.97\times 10^9$, and the value of $c_\std$ for \texttt{3.8.a.a} was roughly $7.08\times 10^{14}$.

\subsection{Results and discussion}  
\begin{figure}
\includegraphics[width=.4\textwidth]{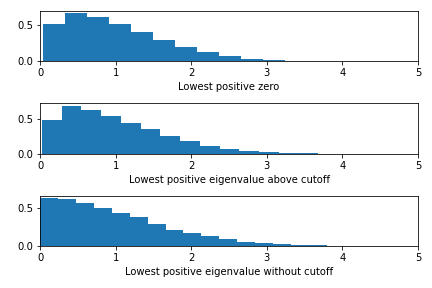}
\includegraphics[width=.4\textwidth]{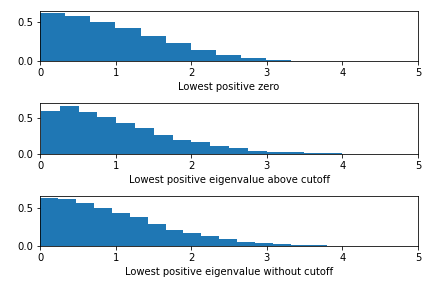}\\

\vspace{.2in}

\includegraphics[width=.4\textwidth]{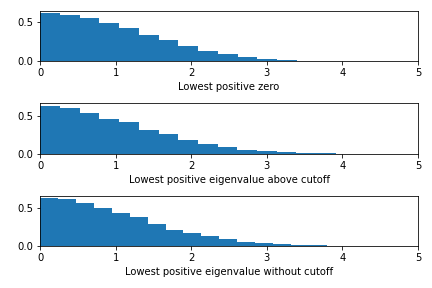}
\includegraphics[width=.4\textwidth]{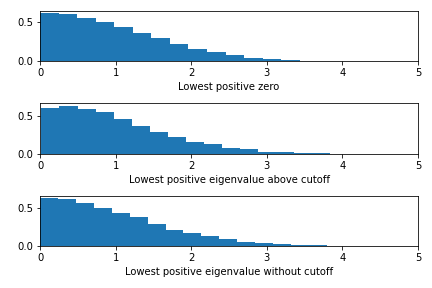}
\caption{Each of these four plots show the distributions of the lowest positive zero (top), the distribution of the lowest eigenvalue for matrices $A$ such that $\abs{\Lambda_A(1,24)}$ is greater than the cutoff we calculated (middle) and the distribution of the lowest positive eigenvalues sampled from the whole of $SO(24)$ (bottom).  From left to right the forms whose central values $f$ are being calculated are \texttt{11.2.a.a}, \texttt{7.4.a.a}, \texttt{3.6.a.a}, and \texttt{3.8.a.a}.  In all cases, we generated random elements from $SO(24)$.}\label{fig:cutoff_std}
\end{figure}

In Figure~\ref{fig:cutoff_std} we see qualitative differences between weight 2 and weights larger than 2.  In particular, in weight 2, we see that the excised model does a better job of modeling the first zero while for weights larger than 2, the full model does better.  On the one hand, this is surprising because it would be expected that a model that takes the arithmetic nature of the modular forms into account would perform better than one that does not.  On the other hand, this is not so surprising because the repulsion that one sees in weight 2 arises from how often the twisted L-functions vanish at the central value.  In weight 4, the excised model overcorrects for something that happens infinitely often but somewhat rarely (according to \cite{CKRS} we expect the number of vanishings for twists up to discriminant $X$ to be on the order of $X^{1/4}\log(X)^{-5/8}$) and in weight 6 and above there is not much difference between the excised model and the full model because the cutoffs are very small (moreover there should not be much repulsion since there are only finitely many vanishings at the central value).

\bibliography{zeros} 

\begin{thebibliography}{10}

\bibitem{bm}
Owen Barrett and Steven~J. Miller.
\newblock An excised orthogonal model for families of cusp forms.
\newblock Preprint.

\bibitem{BaruchMao}
Ehud~Moshe Baruch and Zhengyu Mao.
\newblock {Central value of automorphic L-functions}.
\newblock {\em GAFA Geometric And Functional Analysis}, 17(2):333--384, 2007.

\bibitem{paribook}
Karim Belabas and Henri Cohen.
\newblock {\em Numerical Algorithms for Number Theory: Using Pari/GP}, volume
  254.
\newblock American Mathematical Soc., 2021.

\bibitem{BSP}
Siegfried B{\"o}cherer and Rainer Schulze-Pillot.
\newblock On a theorem of {W}aldspurger and on {E}isenstein series of {K}lingen
  type.
\newblock {\em Math. Ann.}, 288(3):361--388, 1990.

\bibitem{bogomolnykeating1}
Eugene~B Bogomolny and Jon~Peter Keating.
\newblock Random matrix theory and the riemann zeros. i. three-and four-point
  correlations.
\newblock {\em Nonlinearity}, 8(6):1115, 1995.

\bibitem{bogomolnykeating2}
Eugene~B Bogomolny and Jonathan~P Keating.
\newblock Random matrix theory and the riemann zeros ii: n-point correlations.
\newblock {\em Nonlinearity}, 9(4):911, 1996.

\bibitem{code}
Nicol\'{a}s Coloma, Maria Espericueta~Sandoval, Erika Lopez, Francisco Ponce,
  Gustavo Rama, Nathan~C. Ryan, and Alejandro Vargas-Altamirano.
\newblock Repulsion of low-lying zeros of l-functions.
\newblock \url{https://github.com/nathancryan/rmt-families}, 2019.

\bibitem{CKRS}
J.~B. Conrey, J.~P. Keating, M.~O. Rubinstein, and N.~C. Snaith.
\newblock On the frequency of vanishing of quadratic twists of modular
  {$L$}-functions.
\newblock In {\em Number theory for the millennium, {I} ({U}rbana, {IL},
  2000)}, pages 301--315. A K Peters, Natick, MA, 2002.

\bibitem{DHKMS}
E.~Due\~{n}ez, D.~K. Huynh, J.~P. Keating, S.~J. Miller, and N.~C. Snaith.
\newblock A random matrix model for elliptic curve {$L$}-functions of finite
  conductor.
\newblock {\em J. Phys. A}, 45(11):115207, 32, 2012.

\bibitem{Gross}
Benedict~H. Gross.
\newblock Heights and the special values of {$L$}-series.
\newblock In {\em Number theory ({M}ontreal, {Q}ue., 1985)}, volume~7 of {\em
  CMS Conf. Proc.}, pages 115--187. Amer. Math. Soc., Providence, RI, 1987.

\bibitem{ils}
Henryk Iwaniec, Wenzhi Luo, and Peter Sarnak.
\newblock Low lying zeros of families of $ l $-functions.
\newblock {\em Publications Math{\'e}matiques de l'IH{\'E}S}, 91:55--131, 2000.

\bibitem{katzsarnak}
Nicholas Katz and Peter Sarnak.
\newblock Zeroes of zeta functions and symmetry.
\newblock {\em Bulletin of the American Mathematical Society}, 36(1):1--26,
  1999.

\bibitem{KeatingSnaith2}
J.~P. Keating and N.~C. Snaith.
\newblock Random matrix theory and {$L$}-functions at {$s=1/2$}.
\newblock {\em Comm. Math. Phys.}, 214(1):91--110, 2000.

\bibitem{KeatingSnaith1}
J.~P. Keating and N.~C. Snaith.
\newblock Random matrix theory and {$\zeta(1/2+it)$}.
\newblock {\em Comm. Math. Phys.}, 214(1):57--89, 2000.

\bibitem{KZ}
Winfried Kohnen and Don Zagier.
\newblock {Values of $L$-series of modular forms at the center of the critical
  strip}.
\newblock {\em Inventiones mathematicae}, 64(2):175--198, 1981.

\bibitem{MRVT}
Z.~Mao, F.~Rodriguez-Villegas, and G.~Tornar{\'{\i}}a.
\newblock Computation of central value of quadratic twists of modular
  {$L$}-functions.
\newblock In {\em Ranks of elliptic curves and random matrix theory}, volume
  341 of {\em London Math. Soc. Lecture Note Ser.}, pages 273--288. Cambridge
  Univ. Press, Cambridge, 2007.

\bibitem{Mao}
Zhengyu Mao.
\newblock {A generalized Shimura correspondence for newforms}.
\newblock {\em Journal of Number Theory}, 128(1):71--95, 2008.

\bibitem{marshall}
Simon Marshall.
\newblock {Zero repulsion in families of elliptic curve L-functions and an
  observation of Miller}.
\newblock {\em Bulletin of the London Mathematical Society}, 45(1):80--88,
  2013.

\bibitem{mezzadri}
Francesco Mezzadri.
\newblock How to generate random matrices from the classical compact groups.
\newblock {\em Notices Amer. Math. Soc.}, 54(5):592--604, 2007.

\bibitem{miller06}
Steven~J. Miller.
\newblock Investigations of zeros near the central point of elliptic curve
  {$L$}-functions.
\newblock {\em Experiment. Math.}, 15(3):257--279, 2006.
\newblock With an appendix by Eduardo Due\~{n}ez.

\bibitem{montgomery}
Hugh~L Montgomery.
\newblock The pair correlation of zeros of the zeta function.
\newblock In {\em Proc. Symp. Pure Math}, volume~24, pages 181--193, 1973.

\bibitem{odlyzko}
AM~Odlyzko.
\newblock The $10^{22}$-nd zero of the riemann zeta function.
\newblock {\em Dynamical, Spectral, and Arithmetic Zeta Functions: AMS Special
  Session on Dynamical, Spectral, and Arithmetic Zeta Functions, January 15-16,
  1999, San Antonio, Texas}, 290:139, 2001.

\bibitem{os}
Ali~Erhan {\"O}zl{\"u}k and Chip Snyder.
\newblock {On the distribution of the nontrivial zeros of quadratic L-functions
  close to the real axis}.
\newblock {\em Acta Arithmetica}, 91(3):209--228, 1999.

\bibitem{PT2}
Ariel Pacetti and Gonzalo Tornar{\'{\i}}a.
\newblock Examples of the {S}himura correspondence for level {$p^2$} and real
  quadratic twists.
\newblock In {\em Ranks of elliptic curves and random matrix theory}, volume
  341 of {\em London Math. Soc. Lecture Note Ser.}, pages 289--314. Cambridge
  Univ. Press, Cambridge, 2007.

\bibitem{PT1}
Ariel Pacetti and Gonzalo Tornar{\'{\i}}a.
\newblock Shimura correspondence for level {$p^2$} and the central values of
  {$L$}-series.
\newblock {\em J. Number Theory}, 124(2):396--414, 2007.

\bibitem{TR}
Holly Rosson and Gonzalo Tornar{\'{\i}}a.
\newblock Central values of quadratic twists for a modular form of weight 4.
\newblock In {\em Ranks of elliptic curves and random matrix theory}, volume
  341 of {\em London Math. Soc. Lecture Note Ser.}, pages 315--321. Cambridge
  Univ. Press, Cambridge, 2007.

\bibitem{rubinstein}
Michael Rubinstein et~al.
\newblock {Low-lying zeros of L-functions and random matrix theory}.
\newblock {\em Duke Mathematical Journal}, 109(1):147--181, 2001.

\bibitem{lcalc}
Michael~O. Rubinstein.
\newblock {L-calc}.
\newblock \url{https://github.com/nessig/l-calc}, 2019.

\bibitem{rudnicksarnak}
Ze{\'e}v Rudnick, Peter Sarnak, et~al.
\newblock {Zeros of principal L-functions and random matrix theory}.
\newblock {\em Duke Mathematical Journal}, 81(2):269--322, 1996.

\bibitem{pari}
{The PARI~Group}, Univ. Bordeaux.
\newblock {\em {PARI/GP version {\tt 2.11.2}}}, 2018.
\newblock available from \url{http://pari.math.u-bordeaux.fr/}.

\end{thebibliography}
\bibliographystyle{plain}

\end{document}